\documentclass[11pt]{article}%
\usepackage{amsmath}
\usepackage{amsfonts}
\usepackage{amssymb}
\usepackage{graphicx}
\usepackage{a4wide}%
\newtheorem{theorem}{Theorem}

\newtheorem{lemma}[theorem]{Lemma}

\newtheorem{proposition}[theorem]{Proposition}
\newtheorem{remark}[theorem]{Remark}

\begin{document}

\title{On Lagrange multipliers in convex entropy minimization}
\author{Constantin Z\u{a}linescu\thanks{Faculty of Mathematics, University Al.\ I.
Cuza, Bd.\ Carol I, Nr. 11, 700506 Ia\c{s}i, Romania, e-mail:
\texttt{zalinesc@uaic.ro}, and ``O. Mayer'' Institute of Mathematics
of the Romanian Academy, Ia\c{s}i, Romania. }}
\date{}
\maketitle

\begin{abstract}
Based on a characterization of the optimality of a feasible solution of a
convex entropy minimization problem, one shows that the feasible solutions
obtained using formally the Lagrange multipliers method are optimal.

\end{abstract}

\section{\medskip Introduction}

A common procedure to find the solutions of an optimization problem with a
finite number of equality constrains is using the Lagrange multipliers method
(LMM). More precisely, having the function $f:E\subset X\rightarrow\mathbb{R}$
to be minimized (maximized) with the constraints $g_{i}(x)=b_{i},$ where
$g_{i}:E\rightarrow\mathbb{R}$ ($i\in\overline{1,m}$), is to consider the
Lagrangian $L:E\times\mathbb{R}^{m}\rightarrow\mathbb{R}$ defined by
\[
L(x,\lambda):=f(x)+{\sum}_{i=1}^{m}\lambda_{i}(g_{i}(x)-b_{i}),
\]
and to find the critical points $(\overline{x},\overline{\lambda})\in
E\times\mathbb{R}^{m},$ that is
\begin{equation}
\nabla_{x}L(\overline{x},\overline{\lambda})=0,\quad\nabla_{\lambda
}L(\overline{x},\overline{\lambda})=0. \label{r-e0}%
\end{equation}
So, in order to envisage LMM one must have the possibility to speak about
$\nabla_{x}L(\overline{x},\overline{\lambda});$ hence $X$ must be a normed
vector space (or, more generally, a topological vector space), $\overline{x}$
must be in the (algebraic) interior of $E$, and the functions $f$ and $g_{i}$
must be at least G\^{a}teaux differentiable at $\overline{x}.$ Moreover, the
existence of $\overline{\lambda}\in\mathbb{R}^{m}$ verifying the conditions in
Eq.~(\ref{r-e0}) is a necessary condition for the optimality of $\overline{x}$
under supplementary conditions on the data; for a precise statement see for
example \cite[Th.\ 9.3.1]{Lue:69}. Problems appear when the set $E$ has empty
(algebraic) interior, situation in which the differentiability of $f$ and
$g_{i}$ can not be considered (see \cite[pp.~171, 172]{Lue:69}); this is often
the case when $X$ is a function-space, as in entropy minimization (or
maximization) problems. However, in many books and articles on entropy
optimization LMM is used in a formal way. Borwein and Limber (see
\cite{BorLim:92}) describe the main steps of the usual procedure for solving
the entropy minimization problem (see also the survey \cite{Bor:12}); they
mention \textquotedblleft We shall see that this is usually the solution but
each step in the above derivation is suspect and many are wrong without
certain assumptions.'' Pavon and Ferrante (see \cite[Cor.~9.3]{PavFer:13})
establish a sufficient condition for the optimality of the element obtained
using LMM. However, examining their application of this result for
establishing that \textquotedblleft the Gaussian density $p_{c}(x)=(2\pi
)^{-1/2}\exp\big[-\tfrac{1}{2}\tfrac{x^{2}}{\sigma^{2}}\big]$ has maximum
entropy among densities with given mean and variance" we observe that
\cite[Cor.~9.3]{PavFer:13} is not adequate for solving this problem.

The aim of this note is to show that the solutions found using formally LMM
are indeed optimal solutions for the entropy minimization problem%
\[
(EM)\quad\quad\text{minimize }\int_{T}\varphi(x(t))d\mu(t)~~s.t.~~\int_{T}%
\psi_{i}(t)x(t)d\mu(t)=b_{i}~~(i\in\overline{1,m}),
\]
where $\varphi:\mathbb{R}\rightarrow\overline{\mathbb{R}}$ is a proper convex
function, $(T,\mathcal{A},\mu)$ is a measure space, $\psi_{i}:T\rightarrow
\mathbb{R}$ $(i\in\overline{1,m})$ are measurable, and $x\in X$ with $X$ a
linear space of measurable functions. In fact, we provide a characterization
of a solution $\overline{x}$ of $(EM)$ from which we deduce easily that
$\overline{x}$ obtained using LMM is indeed a solution of the problem.

Note that J.M. Borwein and some of his collaborators treated rigorously
problem $(EM)$ when $\mu(T)<\infty$ and the functions $\psi_{i}$ are from
$L_{\infty}(T,\mathcal{A},\mu)$ in a series of papers.

\section{Preliminaries}

Let $(T,\mathcal{A},\mu)$ be a measure space. Set
\[
\mathcal{M}:=\mathcal{M}(T,\mathcal{A},\mu):=\left\{  x:T\rightarrow
\overline{\mathbb{R}}\mid x\text{ is measurable}\right\}  ,
\]
where $\overline{\mathbb{R}}:=\mathbb{R}\cup\{-\infty,\infty\}$ (with
$\infty:=+\infty$), and%
\[
\mathcal{M}_{0}:=\{x\in\mathcal{M}\mid x(t)\in\mathbb{R}\text{ for a.e.\ }t\in
T\},\quad\mathcal{M}_{0}^{+}:=\{x\in\mathcal{M}_{0}\mid x\geq0\text{ a.e.}\}.
\]
As usual we consider as being equal two elements of $\mathcal{M}$ which
coincide almost everywhere (a.e.\ for short). Recall that for every function
$x\in\mathcal{M}$ with values in $\overline{\mathbb{R}}_{+}:=[0,\infty]$ there
exists its integral $\int_{T}xd\mu\in\overline{\mathbb{R}}_{+}$; moreover, if
$\int_{T}xd\mu<\infty,$ then $x\in\mathcal{M}_{0}.$

In the sequel we use the conventions
\[
\infty-\infty:=+\infty+(-\infty):=-\infty+\infty:=\infty,\quad0\cdot(\pm
\infty):=(\pm\infty)\cdot0:=0.
\]
With these conventions $\int_{T}xd\mu:=\int_{T}x_{+}d\mu-\int_{T}x_{-}d\mu$
makes sense for every $x\in\mathcal{M},$ where $\alpha_{+}:=\max\{\alpha,0\}$
and $\alpha_{-}:=(-\alpha)_{+}$ for $\alpha\in\overline{\mathbb{R}}$;
moreover, $\int_{T}xd\mu<\infty$ if and only if $\int_{T}x_{+}d\mu<\infty$ (in
particular $x_{+}\in\mathcal{M}_{0}$), and $\int_{T}xd\mu\in\mathbb{R}$ if and
only if $\int_{T}x_{+}d\mu<\infty$ and $\int_{T}x_{-}d\mu<\infty$ (in
particular $x\in\mathcal{M}_{0}$). The class of those $x\in\mathcal{M}$ with
$\int_{T}xd\mu\in\mathbb{R}$ is denoted, as usual, by $L_{1}(T,\mathcal{A}%
,\mu),$ or simply $L_{1}(T),$ or even $L_{1}.$

\begin{lemma}
\label{lem1} Let $x,y\in\mathcal{M}.$ Then the following assertions hold:

\emph{(a)} If $x\leq y$ then $\int_{T}xd\mu\leq\int_{T}yd\mu.$

\emph{(b)} If $x,y\geq0$ and either $\int_{T}xd\mu<\infty,$ or $\int_{T}%
yd\mu<\infty,$ then $\int_{T}(x-y)d\mu=\int_{T}xd\mu-\int_{T}yd\mu.$

\emph{(c)} If $\int_{T}xd\mu<\infty$ and $\int_{T}yd\mu<\infty$ then $\int%
_{T}(x+y)d\mu=\int_{T}xd\mu+\int_{T}yd\mu$
\end{lemma}

We omit the proof which is standard and uses our conventions.

\begin{remark}
\emph{ Observe that the hypotheses in assertions (b) and (c) of Lemma
\ref{lem1} are essential. For example, taking $T:=\mathbb{R}_{+}$ endowed with
the Lebesgue measure and $x(t):=y(t):=t$ for $t\in T$ we have that $x,y\geq0$
and $\int_{T}xd\mu=\int_{T}yd\mu=\infty,$ while $0=\int_{T}(x-y)d\mu\neq
\int_{T}xd\mu-\int_{T}yd\mu=\infty$; taking $T,$ $\mu,$ $x$ as before and
$z:=-y,$ we have that $\int_{T}xd\mu=\infty,$ $\int_{T}zd\mu=-\infty<\infty$
and $0=\int_{T}(x+z)d\mu\neq\int_{T}xd\mu+\int_{T}zd\mu=\infty.$}
\end{remark}

Consider $\varphi\in\Lambda(\mathbb{R})$ (that is a proper convex function on
$\mathbb{R}$) with $\operatorname*{int}(\operatorname*{dom}\varphi
)\neq\emptyset.$ Because $[\varphi\leq\alpha]:=\{u\in\mathbb{R}\mid
\varphi(u)\leq\alpha\}$ is an interval, it follows immediately that
$\varphi\circ x\in\mathcal{M}$ for every $x\in\mathcal{M},$ where $\varphi
(\pm\infty):=\infty.$ (The notations and notions which are not explained are
standard; see for example \cite{Zal:02}.)

Let us consider a linear space $X\subset\mathcal{M}_{0},$ and define
\begin{equation}
\phi:X\rightarrow\overline{\mathbb{R}},\quad\phi(x):=\int_{T}\varphi\circ
xd\mu. \label{r-phi}%
\end{equation}

\begin{proposition}
\label{p1}Let $\varphi\in\Lambda(\mathbb{R})$ with $\operatorname*{int}%
(\operatorname*{dom}\varphi)\neq\emptyset,$ and let $\phi$ be defined by
(\ref{r-phi}). Then
\[
\operatorname*{dom}\phi=\left\{  x\in X\mid\left(  \varphi\circ x\right)
_{+}\in L_{1}\right\}  \subset\left\{  x\in X\mid x(t)\in\operatorname*{dom}%
\varphi\ \text{a.e.}\right\}
\]
and $\phi$ is convex; in particular, $\operatorname*{dom}\phi$ is convex.
Moreover, if $\varphi$ is strictly convex (on its domain) and $\phi$ is finite
on the convex set $K\subset\operatorname*{dom}\phi$, then $\phi+\iota_{K}$ is
strictly convex, where $\iota_{K}(x):=0$ for $x\in K,$ $\iota_{K}(x):=\infty$
for $x\in X\setminus K.$
\end{proposition}

Proof. The equality follows from our convention $\infty-\infty:=\infty,$ while
the inclusion is obvious. Take $x,y\in\operatorname*{dom}\phi$ and $\lambda
\in{}]0,{}1[{}.$ Since $\varphi$ is convex,
\begin{equation}
\varphi\circ(\lambda x+(1-\lambda)y)\leq\lambda\cdot(\varphi\circ
x)+(1-\lambda)\cdot(\varphi\circ y)\quad\text{a.e.} \label{r1}%
\end{equation}
From Lemma \ref{lem1} (a) we obtain that
\[
\phi(\lambda x+(1-\lambda)y)=\int_{T}\varphi\circ(\lambda x+(1-\lambda
)y)d\mu\leq\int_{T}\left[  \lambda\cdot(\varphi\circ x)+(1-\lambda
)\cdot(\varphi\circ y)\right]  d\mu.
\]
Since $\phi(x)=\int_{T}(\varphi\circ x)d\mu<\infty$ and $\phi(y)=\int%
_{T}(\varphi\circ y)d\mu<\infty,$ using Lemma \ref{lem1} (c) we get
\begin{align*}
\int_{T}\left[  \lambda\cdot(\varphi\circ x)+(1-\lambda)\cdot(\varphi\circ
y)\right]  d\mu &  =\lambda\int_{T}(\varphi\circ x)d\mu+(1-\lambda)\int%
_{T}(\varphi\circ y)d\mu\\
&  =\lambda\phi(x)+(1-\lambda)\phi(y),
\end{align*}
and so $\phi(\lambda x+(1-\lambda)y)\leq\lambda\phi(x)+(1-\lambda)\phi(y).$
Hence $\phi$ is convex.

Let $K$ be a convex subset of $\operatorname*{dom}\phi$ such that
$\phi(K)\subset\mathbb{R}.$ Assume that $\varphi$ is strictly convex, that is
$s,s^{\prime}\in\operatorname*{dom}\varphi$ with $s\neq s^{\prime}$ and
$\lambda\in{}]0,{}1[{}$ imply $\varphi(\lambda s+(1-\lambda)s^{\prime
})<\lambda\varphi(s)+(1-\lambda)\varphi(s^{\prime}).$ Moreover, assume by
contradiction that there exist $x,y\in K$ with $\mu(T_{0})>0,$ where
$T_{0}:=\{t\in T\mid x(t)\neq y(t)\},$ and $\lambda\in{}]0,{}1[{}$ such that
$\phi(\lambda x+(1-\lambda)y)=\lambda\phi(x)+(1-\lambda)\phi(y),$ or,
equivalently $\int_{T}\varphi\circ(\lambda x+(1-\lambda)yd\mu=\int_{T}\left[
\lambda\cdot(\varphi\circ x)+(1-\lambda)\cdot(\varphi\circ y)\right]  d\mu.$
Since $\varphi\circ(\lambda x+(1-\lambda)y)\leq\lambda\cdot(\varphi\circ
x)+(1-\lambda)\cdot(\varphi\circ y)$ a.e.\ and $\varphi\circ(\lambda
x+(1-\lambda)y),$ $\lambda\cdot(\varphi\circ x)$ and $(1-\lambda)\cdot
(\varphi\circ y)$ are from $L_{1},$ it follows that $\varphi\circ(\lambda
x+(1-\lambda)y)=\lambda\cdot(\varphi\circ x)+(1-\lambda)\cdot(\varphi\circ y)$
a.e., contradicting our assumption that $\mu(T_{0})>0.$\hfill$\square$

\medskip

As well known, if $\phi$ takes the value $-\infty,$ then it takes the value
$-\infty$ on on the relative algebraic interior of $\operatorname*{dom}\phi$
denoted $\operatorname*{icr}(\operatorname*{dom}\phi);$ however,
$\operatorname*{icr}(\operatorname*{dom}\phi)$ is empty in many cases of
interest when $X$ is an $L_{p}$ space with $p\in\lbrack1,\infty\lbrack.$

\medskip

Having in view the applications to entropy minimization problems, in the
sequel we consider $\varphi\in\Gamma(\mathbb{R})$ (that is $\varphi\in
\Lambda(\mathbb{R})$ and $\varphi$ is lower semi\-continuous) such that
$\varphi$ is strictly convex on $I:=\operatorname*{dom}\varphi,$
$\operatorname*{int}I\neq\emptyset,$ and $\varphi$ is derivable on
$\operatorname*{int}I;$ this implies that the conjugate $\varphi^{\ast}$ of
$\varphi$ [defined by $\varphi^{\ast}(u)=\sup_{v\in\mathbb{R}}\left(
uv-\varphi(v)\right)  $] is derivable on $\operatorname*{int}%
(\operatorname*{dom}\varphi^{\ast})$ which is nonempty. Moreover, if $a:=\inf
I\in\mathbb{R},$ then either $\varphi(a)=+\infty$ and $\lim_{u\rightarrow
a+}\varphi^{\prime}(u)=-\infty,$ or $\varphi(a)\in\mathbb{R}$ and
$\varphi^{\prime}(a):=\varphi_{+}^{\prime}(a)=\lim_{u\rightarrow a+}%
\varphi^{\prime}(u)$ $(\in\lbrack-\infty,\infty\lbrack).$ Similarly, if
$b:=\sup I\in\mathbb{R},$ then either $\varphi(b)=+\infty$ and $\lim
_{u\rightarrow b-}\varphi^{\prime}(u)=+\infty,$ or $\varphi(b)\in\mathbb{R}$
and $\varphi^{\prime}(b):=\varphi_{-}^{\prime}(b)=\lim_{u\rightarrow
b-}\varphi^{\prime}(u)$ $(\in{}]-\infty,\infty]).$ Assuming that $\phi$ is
proper, then (as seen above) $\phi$ is strictly convex on $\operatorname*{dom}%
\phi.$

\begin{proposition}
\label{p4}Consider $\overline{x},x\in\operatorname*{dom}\phi$ with
$\phi(\overline{x})\in\mathbb{R}$. Then
\begin{equation}
\phi^{\prime}(\overline{x},x-\overline{x}):=\lim_{s\rightarrow0+}\frac
{\phi(\overline{x}+s(x-\overline{x}))-\phi(\overline{x})}{s}=\int_{T}%
\varphi^{\prime}(\overline{x}(t))\cdot\big(x(t)-\overline{x}(t)\big)d\mu(t).
\label{r2}%
\end{equation}

\end{proposition}

Proof. Since $\overline{x},x\in\operatorname*{dom}\phi$ we have that
$\overline{x}(t),x(t)\in\operatorname*{dom}\varphi$ a.e.

Assume first that $\phi(x)\in\mathbb{R}.$ Take $(s_{n})_{n\geq1}\subset{}%
]0,{}1[{}$ a decreasing sequence with $s_{n}\rightarrow0.$ Set
\[
\theta_{n}:=\varphi\circ x-\varphi\circ\overline{x}-\frac{\varphi
\circ(\overline{x}+s_{n}(x-\overline{x}))-\varphi\circ\overline{x}}{s_{n}};
\]
then $0\leq\theta_{n}\leq\theta_{n+1}$ a.e.\ on $T.$ Moreover
\[
\lim_{n\rightarrow\infty}\theta_{n}(t)=\varphi(x(t))-\varphi(\overline
{x}(t))-\varphi^{\prime}(\overline{x}(t))\cdot(x(t)-\overline{x}(t))\in
\lbrack0,+\infty]\text{ for a.e.\ }t\in T.
\]
By Lebesgue's monotone convergence theorem (see \cite[Th.~1.26]{Rud:87}),
$(\varphi^{\prime}\circ\overline{x})\cdot(x-\overline{x})\in\mathcal{M}$ and
\begin{align*}
\phi(x)-\phi(\overline{x})-\phi^{\prime}(\overline{x},x-\overline{x})  &
=\lim_{n\rightarrow\infty}\left[  \phi(x)-\phi(\overline{x})-\frac
{\phi(\overline{x}+s_{n}(x-\overline{x}))-\phi(\overline{x})}{s_{n}}\right]
=\lim_{n\rightarrow\infty}\int_{T}\theta_{n}d\mu\\
&  =\int_{T}\left[  \varphi(x(t))-\varphi(\overline{x}(t))-\varphi^{\prime
}(\overline{x}(t))\cdot(x(t)-\overline{x}(t))\right]  d\mu(t)\in
\overline{\mathbb{R}}_{+}.
\end{align*}
Since $\varphi\circ x,\varphi\circ\overline{x}\in L_{1},$ we get the existence
of $\int_{T}\varphi^{\prime}(\overline{x}(t))\cdot(x(t)-\overline{x}%
(t))d\mu(t)$ as an element of $[-\infty,+\infty\lbrack{},$ and so (\ref{r2}) holds.

Assume now that $\phi(x)=-\infty.$ In this case we have that $(\varphi
^{\prime}\circ\overline{x})\cdot(x-\overline{x})\in\mathcal{M}$, too. With
$(s_{n})_{n\geq1}$ as above, for each $n\geq1$ we have that%
\[
\varphi^{\prime}(\overline{x}(t))\cdot(x(t)-\overline{x}(t))\leq\frac
{\varphi(\overline{x}(t)+s_{n}(x(t)-\overline{x}(t)))-\varphi(\overline
{x}(t))}{s_{n}}\ \ \text{for a.e.}\ \ t\in T. \label{r3}%
\]
Using Lemma \ref{lem1} (a) we get
\[
\int_{T}\varphi^{\prime}(\overline{x}(t))\cdot(x(t)-\overline{x}%
(t))d\mu(t)\leq\frac{\phi(\overline{x}+s_{n}(x-\overline{x}))-\phi
(\overline{x})}{s_{n}}\quad\forall n\geq1,
\]
and so
\begin{equation}
\int_{T}\varphi^{\prime}(\overline{x}(t))\cdot(x(t)-\overline{x}%
(t))d\mu(t)\leq\phi_{+}^{\prime}(\overline{x},x-\overline{x}). \label{r4b}%
\end{equation}
Of course, because $\phi(x)=-\infty$ we have that $\phi\left(  \lambda
x+(1-\lambda)\overline{x}\right)  =\lambda\phi\left(  x)+(1-\lambda
)\phi(\overline{x}\right)  =-\infty$ for $\lambda\in{}]0,{}1[{},$ and so
$\phi_{+}^{\prime}(\overline{x},x-\overline{x})=-\infty.$ From (\ref{r4b}) we
get (\ref{r2}).\hfill$\square$

\section{The entropy minimization problem}

Let us consider $\psi_{1},\ldots,\psi_{m}\in\mathcal{M}_{0}$ and the linear
mappings%
\[
\Psi_{k}:X_{k}\rightarrow\mathbb{R},\quad\Psi_{k}(x):=\int_{T}x\psi_{k}%
d\mu\quad(k\in\overline{1,m}),
\]
where the linear space $X_{k}$ is defined by
\[
X_{k}:=\{x\in\mathcal{M}_{0}\mid x\psi_{k}\in L_{1}\}.
\]
Take also $X_{k}^{0}:=\ker\Psi_{k}:=\{x\in X_{k}\mid\Psi_{k}(x)=0\}$ and set%
\[
\widetilde{X}:=%
{\textstyle\bigcap}
{}_{k=1}^{m}X_{k},\quad\widetilde{X}^{0}:=%
{\textstyle\bigcap}
{}_{k=1}^{m}X_{k}^{0};
\]
note that $\widetilde{X}=\{x\in\mathcal{M}_{0}\mid x\widetilde{\psi}\in
L_{1}\},$ where $\widetilde{\psi}=\left\vert \psi_{1}\right\vert
+\ldots+\left\vert \psi_{m}\right\vert .$

The entropy minimization problem is
\[
(P)\quad\quad\text{minimize }\phi(x)\text{ s.t. }x\in X\cap\widetilde{X}\text{
with }\Psi_{k}(x)=b_{k}~~\forall k\in\overline{1,m},
\]
where $b:=(b_{1},\ldots,b_{m})\in\mathbb{R}^{m}$ is a given element.

Set
\begin{equation}
F_{b}:=\big\{x\in\widetilde{X}\mid\Psi_{k}(x)=b_{k}~\forall k\in\overline
{1,m}\big\}. \label{r-fb}%
\end{equation}
Of course, if $\overline{x}\in F_{b}$ then $F_{b}=\overline{x}+\widetilde{X}%
^{0};$ in particular, $F_{b}$ is a convex set. Because $\varphi$ is strictly
convex, if $\phi+\iota_{F_{b}}$ is proper then $\phi+\iota_{F_{b}}$ is
strictly convex, and so $(P)$ has at most one solution. Said differently, if
$\overline{x}$ is a solution of $(P)$ with $\phi(\overline{x})\in\mathbb{R},$
then $\overline{x}$ is the unique solution of $(P).$

It is known (at least for the Boltzmann--Shannon entropy) that when problem
$(P)$ has a feasible solution $\widetilde{x}\in\operatorname*{dom}\phi$ such
that $\widetilde{x}(t)\in\operatorname*{int}(\operatorname*{dom}\varphi)$ for
a.e.\ $t\in T$ and $\varphi^{\prime}(a)=-\infty$ (when $a=\inf
(\operatorname*{dom}\varphi)\in\mathbb{R}$), $\varphi^{\prime}(b)=+\infty$
(when $b=\sup(\operatorname*{dom}\varphi)\in\mathbb{R}$), if $\overline{x}$ is
the optimal solution of (P) with $\phi(\overline{x})\in\mathbb{R}$, then
$\overline{x}(t)\in\operatorname*{int}(\operatorname*{dom}\varphi)$ for
a.e.\ $t\in T$. Indeed, assume that $a\in\mathbb{R}$ and $\mu(T_{a})>0,$ where
$T_{a}:=\{t\in T\mid\overline{x}(t)=a\}.$ Since $0\leq\phi_{+}^{\prime
}(\overline{x},\widetilde{x}-\overline{x}),$ from (\ref{r2}) we have that
\[
0\leq\int_{T}\varphi^{\prime}(\overline{x}(t))\cdot(\widetilde{x}%
(t)-\overline{x}(t))d\mu(t)\leq\phi(\widetilde{x})-\phi(\overline{x})<\infty,
\]
whence $\int_{T_{a}}\varphi^{\prime}(a)\cdot(\widetilde{x}(t)-a)d\mu
(t)\in\mathbb{R}.$ This is a contradiction, because $\varphi^{\prime
}(a)=-\infty,$ $\widetilde{x}(t)-a>0$ for $t\in T_{a}$ and $\mu(T_{a})>0$. We
get a similar contradiction when $b=\sup(\operatorname*{dom}\varphi
)\in\mathbb{R}$ and $T_{b}:=\{t\in T\mid\overline{x}(t)=b\}$ has positive measure.

\begin{proposition}
\label{p2}Let $\overline{x}\in X\cap F_{b}$ be such that $\phi(\overline
{x})\in\mathbb{R}.$

\emph{(a)} $\overline{x}$ is a solution of problem $(P)$ if and only if (one
of) the following two equivalent conditions hold(s):%
\begin{align}
\phi^{\prime}(\overline{x},x-\overline{x})  &  \geq0\quad\forall x\in
F_{b}\cap\operatorname*{dom}\phi,\label{r50}\\
\int_{T}\varphi^{\prime}(\overline{x}(t))\cdot u(t)d\mu(t)  &  \geq
0\quad\forall u\in K_{\overline{x}}:=\left[  \mathbb{R}_{+}%
(\operatorname*{dom}\phi-\overline{x})\right]  \cap\widetilde{X}^{0}.
\label{r5}%
\end{align}

\emph{(b)} If there exists $\alpha_{1},\ldots,\alpha_{m}\in\mathbb{R}$ such
that $\varphi^{\prime}\circ\overline{x}=\alpha_{1}\psi_{1}+\ldots+\alpha
_{m}\psi_{m},$ then $\overline{x}$ is optimal solution of problem $(P).$
\end{proposition}

Proof. (a) The fact that [$\overline{x}$ is a solution of $(P)$ iff
(\ref{r50}) holds] follows immediately from a known result (see
\cite[Th.\ 3.8]{Jah:07}, \cite[Th.~9.2]{PavFer:13}, \cite[Prop.\ 4]{Zal:89}).
Indeed, from the inequality $\phi^{\prime}(\overline{x},x-\overline{x}%
)\leq\phi(x)-\phi(\overline{x})$ we get the implication $\Leftarrow$. Assume
that $\overline{x}$ is solution of $(P)$ and take $x\in F_{b}\cap
\operatorname*{dom}\phi.$ Then $(1-s)\overline{x}+sx\in F_{b}\cap
\operatorname*{dom}\phi,$ and so $\phi\left(  (1-s)\overline{x}+sx\right)
\geq\phi(\overline{x})$ for $s\in{}]0,{}1[{}.$ Hence $s^{-1}\left[
\phi\left(  (1-s)\overline{x}+sx\right)  -\phi(\overline{x})\right]  \geq0,$
and (\ref{r50}) follows taking the limit for $s\rightarrow0.$

Since $\left(  F_{b}\cap\operatorname*{dom}\phi\right)  -\overline
{x}=(\operatorname*{dom}\phi-\overline{x})\cap\widetilde{X}^{0},$ and using
Proposition \ref{p4}, relation (\ref{r50}) can be rewritten as
\[
\int_{T}\varphi^{\prime}(\overline{x}(t))\cdot u(t)d\mu(t)\geq0\quad\forall
u\in(\operatorname*{dom}\phi-\overline{x})\cap\widetilde{X}^{0},
\]
which, at its turn, is clearly equivalent to (\ref{r5}).

(b) Consider the linear space
\[
Y_{\overline{x}}:=\{u\in\mathcal{M}_{0}\mid(\varphi^{\prime}\circ\overline
{x})\cdot u\in L_{1}\}
\]
and the linear operator
\[
\Theta_{\overline{x}}:Y_{\overline{x}}\rightarrow\mathbb{R},\quad
\Theta_{\overline{x}}(u):=\int_{T}\varphi^{\prime}(\overline{x}(t))\cdot
u(t)d\mu(t).
\]
Since $\phi^{\prime}(\overline{x},x-\overline{x})<\infty$ for every
$x\in\operatorname*{dom}\phi,$ from assertion (a) we have that $\overline{x}$
is a solution of $(P)$ if and only if $K_{\overline{x}}-K_{\overline{x}%
}\subset Y_{\overline{x}}$ and $\Theta_{\overline{x}}(u)\geq0$ for every $u\in
K_{\overline{x}}.$

A sufficient condition for (\ref{r5}) is
\begin{equation}
Y:=X\cap\widetilde{X}\subset Y_{\overline{x}}~~\text{and~~}\Theta
_{\overline{x}}(u)=0~~\forall u\in Y^{0}:=X\cap\widetilde{X}^{0}. \label{r6}%
\end{equation}

Assuming that $Y\subset Y_{\overline{x}},$ condition (\ref{r6}) is equivalent
by \cite[Lem.~3.9]{Rud:91} to the existence of $\alpha_{1},\ldots,\alpha
_{m}\in\mathbb{R}$ such that $\Theta_{\overline{x}}|_{Y}=\alpha_{1}\Psi
_{1}|_{Y}+\ldots+\alpha_{m}\Psi_{m}|_{Y},$ or, equivalently,%
\begin{equation}
\int_{T}\left[  \varphi^{\prime}(\overline{x}(t))-\left(  \alpha_{1}\psi
_{1}(t)+\ldots+\alpha_{m}\psi_{m}(t)\right)  \right]  \cdot u(t)d\mu
(t)=0~~\forall u\in Y. \label{r8}%
\end{equation}
Observing that an obvious sufficient condition for (\ref{r8}) is
\begin{equation}
\varphi^{\prime}(\overline{x}(t))=\alpha_{1}\psi_{1}(t)+\ldots+\alpha_{m}%
\psi_{m}(t)~~\text{for ~a.e.\ }t\in T, \label{r7}%
\end{equation}
the proof is complete.\hfill$\square$

\medskip

It is worth observing that (\ref{r8}) and (\ref{r7}) are equivalent when $\mu$
is $\sigma$-finite and the condition

\medskip

$(H)\quad\forall A\in\mathcal{A}$ with $\mu(A)\in\mathbb{P}:={}]0,\infty
\lbrack{},$ $\exists u\in\mathcal{M}_{0}$ : $u>0$ a.e.\ and $~u\chi_{A}\in Y$

\medskip

\noindent holds. ($\chi_{A}$ is the characteristic function of $A,$ that is
$\chi_{A}(t):=1$ for $t\in A,$ $\chi_{A}(t):=0$ for $t\in T\setminus A$.)
Indeed, the next result holds.

\begin{proposition}
\label{lem2}Let $Y$ verify condition $(H)$ and let $\mu$ be $\sigma$-finite.
Assume that $\overline{y}\in\mathcal{M}$ is such that $\int_{T}\overline
{y}ud\mu=0$ for every $u\in Y.$ Then $\overline{y}=0$ (a.e.).
\end{proposition}

Proof. By contradiction, assume that $\mu([\overline{y}\neq0])>0.$ Setting
$A_{+}:=[\overline{y}>0],$ $A_{-}:=[\overline{y}<0],$ we have that $\mu
(A_{+})>0$ or $\mu(A_{-})>0.$ We may assume that $\mu(A_{+})>0$ (otherwise
replace $\overline{y}$ by $-\overline{y}.$ Because $T$ is $\sigma$-finite,
there exists $A\in\mathcal{A}$ with $A\subset A_{+}$ such that $\mu
(A)\in\mathbb{P}:={}]0,\infty\lbrack{}.$ By our hypothesis, there exists
$u\in\mathcal{M}_{0}$ such that $u>0$ and $u\chi_{A}\in Y.$ It follows that
$\int_{T}\overline{y}\cdot(u\chi_{A})d\mu=0.$ Because $\overline{y}\cdot
(u\chi_{A})\geq0,$ it follows that $\overline{y}\cdot(u\chi_{A})=0$ a.e., and
so $\overline{y}\chi_{A}=0$ a.e.\ This is a contradiction because
$\overline{y}(t)>0$ for every $t\in A$ and $\mu(A)>0.$\hfill$\square$

\medskip

It is worth observing that condition (\ref{r7}) is exactly the one found using
formally LMM. In fact Proposition \ref{p2} and its proof explain how one
arrives rigorously at the sufficient optimality condition of $\overline{x}\in
F_{b}$ with $\phi(\overline{x})\in\mathbb{R}$ in (\ref{r7}). An alternative
justification of this fact in the case of countable sums is done in
\cite{ValZal:16} and applied in \cite{Zal:16}; of course this can also be
obtained using Proposition \ref{p2}~(b) for $T:=\mathbb{N}^{\ast}$ and $\mu$
the counting measure (that is $\mu(A)=\infty$ for $A\subset\mathbb{N}^{\ast}$
infinite and $\mu(A)$ equals the number of elements of $A$ for $A$ finite).

Proposition \ref{p2}~(b) shows that there is no need to verify separately that
the solutions found using LMM in convex or concave entropy optimization are
effectively solutions of $(P)$; this verification is done for example in
\cite[Th.\ 12.1.1]{CovTho:06} and \cite[Ths.\ 3.2, 3.3]{Con}. Note the
following remark from \cite[p.\ 410]{CovTho:06}: \textquotedblleft The
approach using calculus only suggests the form of the density that maximizes
the entropy. To prove that this is indeed the maximum, we can take the second
variation.\textquotedblright

\medskip

As in \cite[Cor.\ 1]{BorChoMar:03}, in the case $\mu(T)<\infty$ and $\psi
_{i}\in L_{\infty}(T)$ $(i\in\overline{1,m}),$ at least for Boltzmann--Shannon
entropy ($\varphi(u):=u\ln u$ for $u\geq0$ with $0\ln0:=0,$ and $\varphi
(u):=\infty$ for $u<0$), for $X=L_{1}(T)$ and $b\in\operatorname*{icr}%
\mathcal{D}$ the problem $(P)$ has optimal solution (provided by LMM), where%
\begin{equation}
\mathcal{D}:=\{b\in\mathbb{R}^{m}\mid F_{b}\cap\operatorname*{dom}\phi
\neq\emptyset\}=\bigg\{b\in\mathbb{R}^{m}\mid\exists x\in\operatorname*{dom}%
\phi,\ \forall i\in\overline{1,m}:\int_{T}x\psi_{i}d\mu=b_{i}\bigg\},
\label{r-D}%
\end{equation}
$F_{b}$ being defined in (\ref{r-fb}).

The situation is completely different in the general case. Let us consider the
problem
\[
(PG)_{m}\quad\quad\text{minimize }\int_{\mathbb{R}}x(t)\ln x(t)dt\text{~ s.t.
}\int_{\mathbb{R}}t^{k-1}x(t)dt=b_{k}~\forall k\in\overline{1,m}%
\]
with $b=(b_{1},\ldots,b_{m})\in\mathbb{R}^{m}$; $(PG)_{m}$ is studied for
example in \cite[Ch.\ 12]{CovTho:06} for $b_{1}=1,$ and in \cite{PavFer:13}
for $m=3$ and $b=(1,0,\sigma^{2})$.

With our previous notation, $T:=\mathbb{R},$ $\mathcal{A}$ is the class of
Lebesgue measurable subsets of $\mathbb{R},$ and $\mu$ is the Lebesgue
measure. Of course, $(PG)_{m}$ is a particular case of problem $(P)$ in which
$\varphi$ is the Boltzmann--Shannon entropy. Of course, $\varphi\in
\Gamma(\mathbb{R}),$ $\operatorname*{dom}\varphi=[0,\infty\lbrack{},$
$\varphi^{\prime}(u)=1+\ln u$ if $u\in\operatorname*{int}(\operatorname*{dom}%
\varphi)={}]0,\infty\lbrack{},$ $\varphi^{\prime}{}(0):=\lim_{u\rightarrow
0+}\varphi^{\prime}(u)=-\infty,$ $\varphi$ is strictly convex on
$\operatorname*{dom}\varphi.$ We take $X:=\mathcal{M}_{0}$ and $\phi$ defined
in (\ref{r-phi}); then $\operatorname*{dom}\phi\subset\mathcal{M}_{0}^{+}.$ In
the present case $\psi_{k}(t)=t^{k-1}$ for $k\in\overline{1,m}$, and so
$X_{k}=\{x\in\mathcal{M}_{0}\mid x\psi_{k}\in L_{1}\};$ hence $X_{1}=L_{1}.$

Let $m=3$; using H\"{o}lder's inequality for $p=q=2$ and $x\psi_{1}$,
$x\psi_{3}$ with $x\in\mathcal{M}_{0}$ we get
\begin{equation}
\int_{\mathbb{R}}\left\vert tx(t)\right\vert dt=\int_{\mathbb{R}}%
\sqrt{\left\vert x(t)\right\vert }\cdot\sqrt{t^{2}\left\vert x(t)\right\vert
}dt\leq\sqrt{\int_{\mathbb{R}}\left\vert x(t)\right\vert dt}\cdot\sqrt
{\int_{\mathbb{R}}t^{2}\left\vert x(t)\right\vert dt}; \label{r12}%
\end{equation}
equality holds in (\ref{r12}) for $x\in\widetilde{X}$ if and only if $x=0$
a.e.\ Hence, if $x\in X_{1}\cap X_{3}$ then $x\in X_{2},$ and so
\[
\widetilde{X}=X_{1}\cap X_{2}\cap X_{3}=X_{1}\cap X_{3}=\{x\in L_{1}\mid
x\psi_{3}\in L_{1}\}.
\]

\begin{proposition}
\label{p-PF}Consider the problem $(PG)_{3}$. Then%
\[
\mathcal{D}=\{(0,0,0)\}\cup\big\{b\in\mathbb{R}^{3}\mid b_{1}>0,\ b_{3}%
>0,\ \left\vert b_{2}\right\vert \leq\sqrt{b_{1}b_{3}}\big\}.
\]
Moreover, if $b=0,$ then $F_{b}\cap\operatorname*{dom}\phi=\{0\},$ and so
$\overline{x}:=0$ is the solution of $(PG)_{3}.$ If $b_{1},b_{3}>0$ and
$\left\vert b_{2}\right\vert <\sqrt{b_{1}b_{3}}$ then the solution and the
value of $(PG)_{3}$ are
\begin{equation}
\overline{x}(t)=\frac{b_{1}^{2}}{\sqrt{2\pi(b_{1}b_{3}-b_{2}^{2})}}%
e^{-\tfrac{1}{2}\frac{b_{1}^{2}}{b_{1}b_{3}-b_{2}^{2}}(t-b_{2})^{2}}%
~~(t\in\mathbb{R}),\quad\phi(\overline{x})=b_{1}\ln\frac{b_{1}^{2}}{\sqrt{2\pi
e(b_{1}b_{3}-b_{2}^{2})}}, \label{r13}%
\end{equation}
respectively. In particular, if $b=(1,0,\sigma^{2})$ with $\sigma>0,$ then
\begin{equation}
\overline{x}(t)=\frac{1}{\sqrt{2\pi\sigma^{2}}}e^{-\frac{t^{2}}{2\sigma^{2}}%
}~~(t\in\mathbb{R}),\quad\phi(\overline{x})=-\ln\sqrt{2\pi e\sigma^{2}}.
\label{r14}%
\end{equation}

\end{proposition}

Proof. Consider $\varphi,$ $X$, $\phi,$ $\psi_{k},$ $X_{k}$ as above. Hence
$\widetilde{X}=X_{1}\cap X_{3}.$

Assume that $F_{b}\cap\operatorname*{dom}\phi\neq\emptyset$ and take $x\in
F_{b}\cap\operatorname*{dom}\phi;$ hence $x\geq0.$ Because $x\psi_{1},$
$x\psi_{3}\geq0,$ it follows that $b_{1},b_{3}\geq0.$ Moreover, from
(\ref{r12}) we obtain that $\left\vert b_{2}\right\vert \leq\sqrt{b_{1}b_{3}%
}.$

If $b_{1}=0,$ then $x\psi_{1}=0$ a.e., and so $x=0$ a.e.; it follows that
$b_{2}=b_{3}=0$ and $F_{b}\cap\operatorname*{dom}\phi=\{0\}.$ The same
conclusion is got when $b_{3}=0.$

Let $b_{1},b_{3}>0$ and assume that $\left\vert b_{2}\right\vert =\sqrt
{b_{1}b_{3}}$. Then equality holds in (\ref{r12}), which implies the existence
of $\alpha\in\mathbb{R}\setminus\{0\}$ such that $x\psi_{1}=\alpha x\psi_{3}.$
Since $\{t\in\mathbb{R}\mid\psi_{1}(t)=\alpha\psi_{3}(t)\}$ is finite, it
follows that $x=0$ a.e., which implies that $b_{1}=0,$ a contradiction.
Therefore, $\left\vert b_{2}\right\vert <\sqrt{b_{1}b_{3}}.$ The fact that
$b_{1},b_{3}>0$ and $\left\vert b_{2}\right\vert <\sqrt{b_{1}b_{3}}$ imply
that $F_{b}\cap\operatorname*{dom}\phi\neq\emptyset$ follows from the fact
that $\overline{x}$ defined in (\ref{r13}) is the optimal solution of
$(PG)_{3},$ as proved below.

Assume that $b_{1},b_{3}>0$ and $\left\vert b_{2}\right\vert <\sqrt{b_{1}%
b_{3}}.$ The problem is to find (if possible) some $\overline{x}\in F_{b}%
\cap\operatorname*{dom}\phi$ $(\subset X\cap\widetilde{X})$ such that
(\ref{r7}) holds. Assuming that such an $\overline{x}$ exists, then
$\overline{x}(t)=e^{c_{0}+c_{1}t+c_{2}t^{2}}$ for $t\in\mathbb{R}$ and some
$c_{0},c_{1},c_{2}\in\mathbb{R}$. Since $\overline{x}\in X_{1}=L_{1},$ we have
necessarily that $c_{2}<0,$ and so $\overline{x}(t)=e^{-\frac{1}{2}%
\alpha(t-\beta)^{2}+\gamma}$ for some $\alpha,\beta,\gamma\in\mathbb{R}$ with
$\alpha>0$ and every $t\in\mathbb{R}$. Imposing $\overline{x}$ to belong to
$F_{b},$ and using the known fact that $\int_{\mathbb{R}}e^{-\tfrac{1}{2}%
t^{2}}dt=\sqrt{2\pi},$ we get
\[
\alpha=\frac{b_{1}^{2}}{b_{1}b_{3}-b_{2}^{2}},\quad\beta=\frac{b_{2}}{b_{1}%
},\quad\gamma=\ln\frac{b_{1}^{2}}{\sqrt{2\pi(b_{1}b_{3}-b_{2}^{2})}}.
\]
Hence $\overline{x}$ is the function defined in (\ref{r13}). Moreover,
\[
\phi(\overline{x})=\int_{\mathbb{R}}\overline{x}(t)\ln\overline{x}%
(t)dt=\int_{\mathbb{R}}\left(  -\frac{1}{2}\alpha(t-\beta)^{2}+\gamma\right)
\overline{x}(t)dt.
\]
Taking into account the constraints and the expressions of $\alpha
,\beta,\gamma$ above, we get the formula for $\phi(\overline{x})$ from
(\ref{r13}).

Moreover, in the general case, for probability densities with mean
$m\in\mathbb{R}$ and variance $\sigma^{2}$ ($\sigma>0$), one has $b_{1}=1,$
$b_{2}=m$ and $b_{3}=\sigma^{2}+2mb_{2}-m^{2}b_{1}=\sigma^{2}+m^{2}.$ From
(\ref{r13}) we get
\[
\overline{x}(t)=\frac{1}{\sqrt{2\pi\sigma^{2}}}e^{-\frac{(t-m)^{2}}%
{2\sigma^{2}}}~~(t\in\mathbb{R}),\quad\phi(\overline{x})=-\ln\sqrt{2\pi
e\sigma^{2}},
\]
which gives (\ref{r14}) when $m=0.$\hfill$\square$

\medskip

As mentioned above, problem $(PG)_{3}$ is considered for $b=(1,0,\sigma^{2})$
in \cite{PavFer:13} and solved applying \cite[Cor.\ 9.3]{PavFer:13}. There
$X=L_{1}(\mathbb{R}),$ whence $X^{\ast}=L_{\infty}(\mathbb{R}),$ and
\[
\mathcal{V}:=\bigg\{x\in X\mid\int_{\mathbb{R}}x(t)dt=\int_{\mathbb{R}%
}tx(t)dt=\int_{\mathbb{R}}t^{2}x(t)dt=0\bigg\}.
\]
It is not explained why \cite[(9.6)]{PavFer:13} holds and which is the
annihilator of $\mathcal{V}$ in order to take $\overline{x}$ of the form
$t\mapsto Ce^{\vartheta_{1}t+\vartheta_{2}t^{2}}.$

Proposition \ref{p-PF} provides an example in which $\mu(T)=\infty$ and the
problem $(P)$ has optimal solutions for all $b\in\mathcal{D}.$ In
\cite{BorChoMar:03} it is presented a situation with $X=L_{1}(0,\infty)$ and
$\varphi$ the Boltzmann--Shannon entropy in which $(P)$ has optimal solutions
for all $b\in\operatorname*{icr}\mathcal{D},$ as in the case $\mu(T)<\infty$
and $\psi_{i}\in L_{\infty}(T).$

Problem $(P)$ is considered in \cite{ValZal:16} and \cite{Zal:16} for
$T:=\mathbb{N}^{\ast}$, $\mu$ the counting measure, and $\varphi(u)=u\ln u-u$
for $u\geq0,$ $\varphi(u)=\infty$ for $u<0.$ Practically, a complete study of
$(P)$ for $m=1$ is given in \cite[Prop.\ 3.3]{ValZal:16}; so, besides
providing the value of problem $(P)$ for $b\in\mathcal{D}$ [$\mathcal{D}$
being defined in (\ref{r-D})], when $\mathcal{D}\neq\{0\}$ it is shown that
either $(P)$ has optimal solution for each $b\in\mathcal{D},$ or
$(\operatorname*{int}\mathcal{D})\setminus\{b\in\mathcal{D}\mid(P)$ has
optimal solution$\}$ is nonempty. In \cite[Prop.\ 3.4]{ValZal:16}, for $m=2$
one has an example in which $\mathcal{D}=\{(0,0)\}\cap\big(  (0,\infty
)\times\mathbb{R}\big) $ and for every $b_{1}>0$ there exists only one
$b_{2}\in\mathbb{R}$ for which $(P)$ has optimal solution which (moreover) can
be found using formally LMM. A complete solution of problem $(P)$ for $m=2$
and $\psi_{1}\equiv1$ is given in \cite[Th.\ 4.1]{Zal:16}; the conclusions are
similar to those in \cite[Prop.\ 3.3]{ValZal:16} presented above.

The study of problem $(P)$ for arbitrary measure spaces is done by P.
Mar\'{e}chal in \cite{Mar:99, Mar:01} using a duality approach in which the
primal space is similar to $\widetilde{X}.$


\begin{thebibliography}{99}                                                                                               %


\bibitem {Bor:12}J.M. Borwein, \emph{Maximum entropy and feasibility methods
for convex and nonconvex inverse problems}, Optimization 61 (2012), 1--33.

\bibitem {BorChoMar:03}J.M. Borwein, R. Choksi, P. Marechal, \emph{Probability
distributions of assets inferred from option prices via the principle of
maximum entropy}, SIAM J. Optim. 14 (2003), 464--478

\bibitem {BorLim:92}J. M. Borwein, M. A. Limber, \emph{On entropy maximization
via convex programming}, Preprint CORR 92-16, University of Waterloo (1992).

\bibitem {Con}K. Conrad, \emph{Probability distributions and maximum entropy},
Expository paper, www.math.uconn.edu/\symbol{126}kconrad/blurbs/analysis/entropypost.pdf

\bibitem {CovTho:06}T. M. Cover, J. A. Thomas, Elements of Information Theory,
2nd ed., John Wiley \& Sons, Hoboken, NJ, 2006.

\bibitem {Jah:07}J. Jahn, \emph{Introduction to the Theory of Nonlinear
Optimization}, 3rd ed., Springer, Berlin, 2007.

\bibitem {Lue:69}D. Luenberger, \emph{Optimization by Vector Space Methods},
JohnWiley \& Sons, Inc. (1969).

\bibitem {Mar:99}P. Marechal, \emph{On the principle of maximum entropy on the
mean as a methodology for the regularization of inverse problems},
Grigelionis, B. (ed.) et al., Probability theory and mathematical statistics.
Proceedings of the 7th international Vilnius conference, Vilnius, Lithuania,
August, 12-18, 1998. Vilnius: TEV. 481-492 (1999).

\bibitem {Mar:01}P. Marechal, \emph{A note on entropy optimization}, in
``Approximation, Optimization and Mathematical Economics'', M. Lassonde ed.,
Physica-Verlag, Heidelberg, 2001, pp.\ 205--211.

\bibitem {PavFer:13}M. Pavon, A. Ferrante, \emph{On the geometry of maximum
entropy problems}, SIAM Review 55 (2013), 415--439.

\bibitem {Rud:87}W. Rudin, \emph{Real and Complex Analysis}, (3rd edition),
McGraw-Hill, Inc., 1987.

\bibitem {Rud:91}W. Rudin, \emph{Functional Analysis} (2nd edition),
McGraw-Hill, Inc., 1991.

\bibitem {ValZal:16}C. Vall\'{e}e, C. Z\u{a}linescu, \emph{Series of convex
functions: subdifferential, conjugate and applications to entropy
minimization}, J. Convex Anal. 2016.

\bibitem {Zal:89}C. Z\u{a}linescu, \emph{On Gwinner's paper \textquotedblleft
Results of Farkas type\textquotedblright}, Numer. Funct. Anal. Optim. 10
(1989), 199--210.

\bibitem {Zal:02}C. Z\u{a}linescu, \emph{Convex Analysis in General Vector
Spaces}, World Scientific, River Edge (NJ), 2002.

\bibitem {Zal:16}C. Z\u{a}linescu, \emph{On the entropy minimization problem
in Statistical Mechanics}, J. Math. Anal. Appl. 2016. DOI: 10.1016/j.jmaa.2016.10.020
\end{thebibliography}
\end{document}